\documentclass{article}

\usepackage[T1]{fontenc}
\usepackage[utf8]{inputenc}

\usepackage{lmodern}
\usepackage{microtype}
\usepackage{amsmath,amssymb,mathtools,amsthm}
\usepackage{booktabs}
\usepackage{graphicx}
\usepackage{float}
\usepackage{placeins}
\usepackage[hidelinks,pdfpagemode=UseNone,bookmarksopen=false]{hyperref}

\usepackage[margin=1in]{geometry}
\allowdisplaybreaks
\newtheorem{theorem}{Theorem}[section]

\usepackage{cite}
\newcommand{\Gtilde}{\widetilde{G}}
\newcommand{\Ntilde}{\widetilde{N}}
\newcommand{\Dop}{\mathcal{D}}
\newcommand{\tailF}{\overline{F}}

\title{Efficient Recursive Evaluation of Gaussian, Gamma, and Inverse Gaussian Tail Probabilities}
\author{Joshua Joseph George and Hassan Safouhi\footnote{Corresponding author: hsafouhi@ualberta.ca.}\\
Campus Saint-Jean, University of Alberta\\
8406 91 Street, Edmonton (AB) T6C 4G9, Canada}
\date{}

\begin{document}
\maketitle

\vspace*{0.5cm}

\vspace*{0.5cm}
{\bf \large Abstract}. \hskip 0.10cm Accurate computation of tail probabilities motivates efficient high-order implementations of the $G_n^{(1)}$ extrapolation transformation. We derive finite recurrence relations in the transformation order for Gaussian, gamma, and inverse Gaussian tail approximations. The Gaussian and inverse Gaussian numerator and denominator sequences satisfy four-term recurrences, whereas the gamma sequences satisfy a three-term recurrence. The derivation exploits the structure of the underlying differential operator and shows that the reduced numerator and denominator obey recurrences with identical coefficients. This construction avoids the repeated evaluation of higher-order derivatives and nested sums required by the explicit formulations. As a result, generating the approximation sequence through order $n$ requires $O(n)$ arithmetic operations, compared with $O(n^3)$ for the Gaussian and gamma distributions and $O(n^4)$ for the inverse Gaussian distribution. The recurrence algorithms also require only a fixed number of previously computed quantities at each order. Numerical results for the three distributions are in excellent agreement with values reported in the literature.

\vspace*{0.5cm}

\section{Introduction}
For a probability density function $f$ supported on $(a,\infty)$, where $a$ may be $-\infty$, the upper-tail probability is:
\begin{equation}
\tailF(x)=\int_x^\infty f(t)\,dt.
\label{tail}
\end{equation}
Accurate numerical evaluation of \eqref{tail} is important in statistical inference, reliability analysis, and first-passage modelling, especially when the same distribution must be evaluated repeatedly at different arguments \cite{Good-23-243-86,Gray-Wang-86-159-91,Chhikara-Folks-89,Ryan-07}. For small upper-tail probabilities, the accuracy of the evaluation and its arithmetic cost both matter. Gaussian, gamma, and inverse Gaussian tails provide three representative classes of integrals with distinct exponential and algebraic factors.

The $G$ transformation was introduced in \cite{Gray-Atchison-4-363-67} and developed in \cite{Gray-tchison-McWilliams-8-365-71,Gray-Wang-29-271-92}. The Slevinsky--Safouhi formulae for higher-order differentiation \cite{Safouhi37} facilitated the algorithm developed in \cite{Safouhi40} for the $G_n^{(1)}$ transformation through successive applications of a differential operator. Although recursive in its mathematical formulation, that approach requires the evaluation of higher-order derivatives, leading to increasing computational cost. It was subsequently applied to the accurate computation of tail probabilities for several distributions \cite{Safouhi43}. A further advance was achieved for incomplete Bessel functions through finite recurrence relations \cite{Safouhi62}, which avoid repeated higher-order differentiation and reduce the cost of generating their approximation sequence from $O(n^4)$ to $O(n)$. The present work builds upon this development by deriving explicit distribution-specific recurrence formulas for Gaussian, gamma, and inverse Gaussian tail probabilities.

The present contribution consists of explicit finite recurrences for the denominator and reduced-numerator sequences of the Gaussian, gamma, and inverse Gaussian $G_n^{(1)}$ tail approximations. Both sequences satisfy four-term recurrences for the Gaussian and inverse Gaussian distributions, whereas three-term recurrences suffice for the gamma distribution. A central feature of the derivation is the exact cancellation of the additional term arising in the inhomogeneous differential propagation of the reduced numerator. Consequently, the reduced numerator and denominator satisfy recurrences with identical coefficients. Each new order can thus be obtained from a fixed number of preceding quantities without re-evaluating the nested sums.

The numerical performance of the proposed recurrences is examined using the 35 test cases considered in \cite{Safouhi43}. Following the same approach, convergence is assessed through the relative difference between successive $G_n^{(1)}$ approximations, with a prescribed tolerance of $10^{-15}$. The computed results show excellent agreement with those reported in \cite{Safouhi43}. The numerical results demonstrate the convergence and computational efficiency of the proposed recurrence algorithms across a range of distribution parameters. The recurrences preserve the approximants of the original transformation while reducing computational complexity and requiring only a fixed number of working quantities.

\section{Recurrence formulation and derivation}
We recall only the essential definitions and equations required for the derivation of the finite recurrences. For a detailed development of the $G_n^{(1)}$ construction and its application to tail probabilities, we refer the reader to \cite{Safouhi40,Safouhi43}.

We let:
\begin{equation}
F(x)=\int_a^x f(t)\,dt
\quad \text{and}\quad
I=\int_a^\infty f(t)\,dt
\quad \Longrightarrow\quad
\tailF(x)=I-F(x).
\label{Ftail}
\end{equation}

For the $m=1$ case of the $G$ transformation, we set:
\begin{equation}
\Dop=x^2\frac{d}{dx}, \qquad
D_0(x)=\frac{1}{x^{\sigma_0}f(x)}
\qquad \text{and}\qquad
N_0(x)=F(x)D_0(x),
\label{initial}
\end{equation}
and we define for $n \geq 1$:
\begin{equation}
D_n = \Dop D_{n-1}
\quad \text{and}\quad
N_n = \Dop N_{n-1}.
\label{ND1}
\end{equation}

As given in \cite{Safouhi43}, the $G_n^{(1)}$ approximation of the semi-infinite integral $I$ is:
\begin{equation}
G_n^{(1)}=\frac{N_n}{D_n}.
\label{ND}
\end{equation}

For the three distributions considered below, the $m=1$ construction gives $\sigma_0=-1$ for the standardized Gaussian density and $\sigma_0=0$ for the gamma and inverse Gaussian densities~\cite{Safouhi43}.

The reduced numerator and the corresponding upper-tail approximation $\Gtilde_n^{(1)}$ of $\tailF(x)$ are:
\begin{equation}
\Ntilde_n=N_n-FD_n
\quad \text{and}\quad
\Gtilde_n^{(1)}=\frac{\Ntilde_n}{D_n}.
\label{reduced}
\end{equation}

Since $F'(x)=f(x)$, the two sequences satisfy:
\begin{equation}
D_{n+1}=\Dop D_n
\quad \text{and}\quad
\Ntilde_{n+1}=\Dop\Ntilde_n+x^2 f(x) D_n.
\label{onestep}
\end{equation}

Thus $\Ntilde_n$ is not generated by the homogeneous iteration satisfied by $D_n$.

Suppose that, for a given distribution, the denominator satisfies:
\begin{equation}
D_{n+1}=A_nD_n+B_nD_{n-1}+C_nD_{n-2},
\label{generic}
\end{equation}
where $C_n$ may vanish. Applying $\Dop=x^2\dfrac{d}{dx},$ and using \eqref{onestep} gives:
\begin{equation}
\begin{aligned}
D_{n+2}={}&A_nD_{n+1}+(\Dop A_n+B_n)D_n +(\Dop B_n+C_n)D_{n-1}+(\Dop C_n)D_{n-2}.
\end{aligned}
\label{genericdiff}
\end{equation}

The recurrence holds at the next order if, for some $h(x)$:
\begin{equation}
\begin{aligned}
A_{n+1}&=A_n+h(x),\\
\Dop A_n+B_n&=h(x) A_n+B_{n+1},\\
\Dop B_n+C_n&=h(x) B_n+C_{n+1},\\
\Dop C_n&=h(x) C_n.
\end{aligned}
\label{closure}
\end{equation}

For the reduced numerator, we write \eqref{onestep} as:
\begin{equation}
\Dop\Ntilde_j=\Ntilde_{j+1}- x^2f(x) D_j.
\label{forced}
\end{equation}

We assume at index $n$ that:
\begin{equation}
\Ntilde_{n+1}=A_n\Ntilde_n+B_n\Ntilde_{n-1}+C_n\Ntilde_{n-2}.
\label{numinduction}
\end{equation}

Applying $\Dop$ to \eqref{numinduction}, using the product rule, and substituting \eqref{forced} for each differentiated reduced numerator gives:
\begin{equation}
\begin{aligned}
\Dop\Ntilde_{n+1} & = A_n\Ntilde_{n+1}+(\Dop A_n+B_n)\Ntilde_n +(\Dop B_n+C_n)\Ntilde_{n-1}+(\Dop C_n)\Ntilde_{n-2}\\
& - x^2f(x) \left(A_nD_n+B_nD_{n-1}+C_nD_{n-2}\right).
\end{aligned}
\label{numdiff}
\end{equation}

Using the denominator recurrence, the last term becomes:
\begin{equation}
- x^2f(x) \left(A_nD_n+B_nD_{n-1}+C_nD_{n-2}\right)=- x^2f(x) D_{n+1},
\label{cancellation1}
\end{equation}
whereas the left-hand side of \eqref{numdiff} is:
\begin{equation}
\Dop\Ntilde_{n+1}=\Ntilde_{n+2}- x^2f(x) D_{n+1}.
\label{cancellationlhs}
\end{equation}

The contribution of the additional term therefore cancels exactly. After using \eqref{closure}, the remaining terms give:
\begin{equation}
\begin{aligned}
\Ntilde_{n+2}={}&A_n\Ntilde_{n+1}+(h(x) A_n+B_{n+1})\Ntilde_n +(h(x) B_n+C_{n+1})\Ntilde_{n-1}+ h(x) C_n\Ntilde_{n-2}\\
={}&A_{n+1}\Ntilde_{n+1}+B_{n+1}\Ntilde_n+C_{n+1}\Ntilde_{n-1}.
\end{aligned}
\label{cancellation2}
\end{equation}

It follows that $D_n$ and $\Ntilde_n$ satisfy the same finite recurrence, provided that the recurrence is verified at its initial index. The ratio in \eqref{reduced} is defined whenever $D_n\ne0$.

We now apply the recurrence to the Gaussian, gamma, and inverse Gaussian distributions.

\subsection{Gaussian recurrence}
For the standard normal density and $x>0$:
\begin{equation}
\phi(x)=\frac{1}{\sqrt{2\pi}}e^{-x^2/2},
\label{normalinitial}
\end{equation}
and:
\begin{equation}
\left\{
\begin{array}{lll}
D_0 = x \sqrt{2\pi}e^{x^2/2} & \text{and} & D_1 =(x^3+x)D_0\\[1ex]
\Ntilde_0=0 & \text{and} & \Ntilde_1=x^3.
\end{array}
\right.
\label{normalinitialPlus}
\end{equation}

For a nonstandard normal variable, we use the standardized argument $z=(x-\mu)/\sigma>0$ in the recurrence, with $\overline F(\mu)=1/2$ and the identity $\overline F(x)=1-\overline F(2\mu-x)$ when $z<0$. In the Gaussian formulas below, $x$ denotes this positive standardized argument.

\begin{theorem}[Gaussian recurrence]
\label{thm:normal}
We set $Q_{-1}=0$. For either $Q_n=D_n$ or $Q_n=\Ntilde_n$, and for $n\ge1$:
\begin{equation}
\begin{aligned}
Q_{n+1}={}&\left[x^3+(3n+1)x\right]Q_n
-n(3n-1)x^2Q_{n-1}+n(n-1)^2x^3Q_{n-2}.
\end{aligned}
\label{normalRec}
\end{equation}
\end{theorem}
\begin{proof}

We set:
\begin{equation}
A_n=x^3+(3n+1)x, \quad
B_n=-n(3n-1)x^2
\quad \text{and}\quad
C_n=n(n-1)^2x^3.
\label{normalABC}
\end{equation}

The shift of the leading coefficient is:
\begin{equation}
A_{n+1}-A_n=3x \Longrightarrow \quad h=3x.
\label{normalshift}
\end{equation}

Moreover:
\begin{equation}
\Dop A_n=3x^4+(3n+1)x^2
\quad \text{and}\quad
\Dop B_n=-2n(3n-1)x^3,
\label{normaldiffcoeff}
\end{equation}

For the first closure identity, the two sides reduce to the same polynomial:
\begin{equation}
\begin{aligned}
\Dop A_n+B_n
&=3x^4+(-3n^2+4n+1)x^2,\\
3xA_n+B_{n+1}
&=3x^4+(-3n^2+4n+1)x^2.
\end{aligned}
\label{normalclosurecheck1}
\end{equation}

Similarly:
\begin{equation}
\Dop B_n+C_n=n(n^2-8n+3)x^3
=3xB_n+C_{n+1}.
\label{normalclosurecheck2}
\end{equation}

By collecting like powers of $x$, we obtain:
\begin{equation}
\Dop A_n+B_n=3xA_n+B_{n+1}
\quad \text{and}\quad
\Dop B_n+C_n=3xB_n+C_{n+1}.
\label{normalclosure12}
\end{equation}

The last closure condition is also satisfied:
\begin{equation}
\Dop C_n=3n(n-1)^2x^4=3xC_n.
\label{normalclosure3}
\end{equation}

It follows that \eqref{closure} holds with $h=3x$. For the base index, differentiation of $D_1=(x^3+x)D_0$ gives:
\begin{equation}
D_2=(x^3+4x)D_1-2x^2D_0.
\label{normalbaseD}
\end{equation}

Similarly, using $\Ntilde_1=x^3$ in \eqref{onestep} gives:
\begin{equation}
\Ntilde_2=x^6+4x^4=(x^3+4x)\Ntilde_1.
\label{normalbaseN}
\end{equation}

These are exactly the $n=1$ cases of \eqref{normalRec}. The cancellation in \eqref{cancellation2} then proves the recurrence for all $n\geq 1$.
\end{proof}

\subsection{Gamma recurrence}
For the gamma density with shape $a>0$ and scale $b>0$:
\begin{equation}
f_\Gamma(x)=\frac{x^{a-1}e^{-x/b}}{\Gamma(a)b^a},
\label{gammadensity}
\end{equation}
and:
\begin{equation}
\left\{
\begin{array}{lll}
D_0=\Gamma(a)b^a x^{1-a}e^{x/b} & \text{and} &D_1=\left[(1-a)x+\frac{x^2}{b}\right]D_0\\[1ex]
\Ntilde_0=0 & \text{and} &\Ntilde_1=x^2.
\end{array}
\right.
\label{gammainitial}
\end{equation}

\begin{theorem}[Gamma recurrence]
\label{thm:gamma}
For $n\ge1$ and either $Q_n=D_n$ or $Q_n=\Ntilde_n$:
\begin{equation}
Q_{n+1}=\left[(2n+1-a)x+\frac{x^2}{b}\right]Q_n
-n(n-a)x^2Q_{n-1}.
\label{gammaRec}
\end{equation}
\end{theorem}
\begin{proof}
With:
\begin{equation}
A_n=(2n+1-a)x+\frac{x^2}{b}
\quad \text{and}\quad
B_n=-n(n-a)x^2,
\label{gammaAB}
\end{equation}
the leading coefficient satisfies:
\begin{equation}
A_{n+1}-A_n=2x,
\label{gammashift}
\end{equation}
and differentiation gives:
\begin{equation}
\Dop A_n=(2n+1-a)x^2+\frac{2x^3}{b}
\quad \text{and}\quad
\Dop B_n=-2n(n-a)x^3.
\label{gammadiffcoeff}
\end{equation}

The coefficient of $x^2$ in $\Dop A_n+B_n$ is:
\begin{equation}
(2n+1-a)-n(n-a)=-n^2+(a+2)n+1-a,
\label{gammacoeffcheck}
\end{equation}
which is also the coefficient obtained from $2xA_n+B_{n+1}$.

Consequently:
\begin{equation}
\Dop A_n+B_n=2xA_n+B_{n+1}
\quad \text{and}\quad
\Dop B_n=2xB_n.
\label{gammaclosure}
\end{equation}

It follows that \eqref{closure} holds with $h=2x$ and $C_n=0$. Therefore, the gamma recurrence contains only three terms. The denominator identity at $n=1$ follows by differentiating $D_1$ in \eqref{gammainitial}. For the reduced numerator, \eqref{onestep} gives:
\begin{equation}
\Ntilde_2=2x^3+x^2\left[(1-a)x+\frac{x^2}{b}\right]
=\left[(3-a)x+\frac{x^2}{b}\right]\Ntilde_1.
\label{gammaBaseN}
\end{equation}

Thus, the recurrence \eqref{gammaRec} holds at $n=1$ for both $D_n$ and $\widetilde N_n$. Equation \eqref{cancellation2} establishes the induction step, proving that the recurrence holds for all $n \geq 1$.
\end{proof}

\subsection{Inverse Gaussian recurrence}
For the inverse Gaussian distribution with $\mu>0$ and $\lambda>0$:
\begin{equation}
f_{\mathrm{IG}}(x)=\left(\frac{\lambda}{2\pi x^3}\right)^{1/2}
\exp\left[-\frac{\lambda(x-\mu)^2}{2\mu^2x}\right]
\quad \text{for } \quad x>0.
\label{igdensity}
\end{equation}
The initial values are:
\begin{equation}
\left\{
\begin{array}{lll}
D_0=1/f_{\mathrm{IG}} & \text{and} & D_1 = \displaystyle \left[\frac{\lambda x^2}{2\mu^2} + \frac32x-\frac{\lambda}{2}\right] D_0\\[2ex]
\Ntilde_0=0 & \text{and} &\Ntilde_1=x^2.
\end{array}
\right.
\label{iginitial}
\end{equation}

\begin{theorem}[Inverse Gaussian recurrence]
\label{thm:ig}
We set $Q_{-1}=0$. For $n\ge1$ and either $Q_n=D_n$ or $Q_n=\Ntilde_n$:
\begin{equation}
\begin{aligned}
Q_{n+1} & = \left[\frac{\lambda x^2}{2\mu^2}+\left(\frac32+2n\right)x-\frac{\lambda}{2}\right]Q_n +\left[\lambda nx-n(n-1)x^2-\frac32nx^2\right]Q_{n-1} -\frac{\lambda}{2}n(n-1)x^2Q_{n-2}.
\end{aligned}
\label{igRec}
\end{equation}
\end{theorem}

\begin{proof}
We set:
\begin{equation}
\begin{aligned}
A_n &=\frac{\lambda x^2}{2\mu^2}+\left(\frac32+2n\right)x-\frac{\lambda}{2},\\
B_n &=\lambda nx-n(n-1)x^2-\frac32nx^2,\\
C_n &=-\frac{\lambda}{2}n(n-1)x^2.
\end{aligned}
\label{igABC}
\end{equation}

The leading shift is:
\begin{equation}
A_{n+1}-A_n=2x,
\label{igshift}
\end{equation}
while:
\begin{equation}
\Dop A_n=\frac{\lambda x^3}{\mu^2}+\left(\frac32+2n\right)x^2,
\label{igdiffA}
\end{equation}
and:
\begin{equation}
\Dop B_n=\lambda n x^2-2n(n-1)x^3-3nx^3.
\label{igdiffB}
\end{equation}

Adding the recurrence coefficients gives:
\begin{equation}
\begin{aligned}
\Dop A_n+B_n
={}&\frac{\lambda x^3}{\mu^2}+\lambda nx +\left(\frac32+2n-n(n-1)-\frac32n\right)x^2,
\end{aligned}
\label{igclosurecheckA}
\end{equation}
while:
\begin{equation}
\begin{aligned}
\Dop B_n+C_n
={}&\left(\lambda n-\frac{\lambda}{2}n(n-1)\right)x^2 -\left[2n(n-1)+3n\right]x^3.
\end{aligned}
\label{igclosurecheckB}
\end{equation}

Expanding $2xA_n+B_{n+1}$ and $2xB_n+C_{n+1}$ gives the same coefficients. Therefore:
\begin{equation}
\Dop A_n+B_n=2xA_n+B_{n+1}
\quad \text{and}\quad
\Dop B_n+C_n=2xB_n+C_{n+1}.
\label{igclosure12}
\end{equation}

Finally:
\begin{equation}
\Dop C_n=-\lambda n(n-1)x^3=2xC_n.
\label{igclosure3}
\end{equation}

It follows that \eqref{closure} holds with $h=2x$. The denominator base relation follows by differentiating $D_1$ in \eqref{iginitial}. For the reduced numerator, \eqref{onestep} gives:
\begin{equation}
\Ntilde_2=2x^3+x^2\left[\frac{\lambda x^2}{2\mu^2}+\frac32x-\frac\lambda2\right]
=\left[\frac{\lambda x^2}{2\mu^2}+\frac72x-\frac\lambda2\right]\Ntilde_1.
\label{igBaseN}
\end{equation}
This is the $n=1$ instance of \eqref{igRec}, and the cancellation of the additional term proves the recurrence for all $n\geq1$.
\end{proof}

\section{Computational aspects and complexity}
\label{sec_complexity}

For a given distribution and evaluation point, we first compute $D_0$ and $D_1$, set $\Ntilde_0=0$ and compute the corresponding $\Ntilde_1$, and then propagate both sequences using the same recurrence coefficients. At each order:
\begin{equation}
\Gtilde_k^{(1)}=\frac{\Ntilde_k}{D_k}.
\label{implementratio}
\end{equation}

The gamma recurrence requires only the two preceding states, whereas the Gaussian and inverse Gaussian recurrences require three. For adaptive calculations, successive approximations can be compared to assess convergence.

The numerator and denominator may become large even when their ratio remains moderate. To avoid forming a potentially large $D_0$, one can propagate the normalized sequences:
\begin{equation}
 d_n=\frac{D_n}{D_0},\qquad u_n=\frac{\Ntilde_n}{D_0},\qquad
 d_0=1,\quad d_1=\frac{D_1}{D_0},\quad
 u_0=0,\quad u_1=x^2 f(x),\quad
 \Gtilde_n^{(1)}=\frac{u_n}{d_n}.
\label{normalizedStates}
\end{equation}
Both normalized sequences obey the corresponding recurrence with the same coefficients. In the Gaussian case, $x$ and $f$ in \eqref{normalizedStates} denote the positive standardized argument and standard normal density. 

The recurrence algorithm requires $O(1)$ working storage, whereas storing the complete sequence of approximants requires $O(n)$ storage.

The arithmetic-operation count follows directly from the finite recurrences and refers to fixed-precision scalar operations. It compares the recurrences with direct evaluations of the earlier nested-sum representations. For the Gaussian and gamma distributions, the explicit reduced numerator at fixed order $k$ contains an outer sum over $r=1,\ldots,k$ and coefficient sums with lengths $O(k-r)$ and $O(r)$.

For these direct nested-sum evaluations, the number of scalar summands at one order is:
\begin{equation}
\sum_{r=1}^{k}O\big((k-r)+r\big)=O(k^2).
\label{oldfixedNG}
\end{equation}

Generating all orders through $n$ therefore requires:
\begin{equation}
T_{\mathrm{exp}}^{N,\Gamma}(n)=\sum_{k=1}^{n}O(k^2)=O(n^3).
\label{oldNG}
\end{equation}

The formula for the inverse Gaussian contains one additional nesting level, so one order requires $O(k^3)$ work and the complete sequence requires:
\begin{equation}
T_{\mathrm{exp}}^{\mathrm{IG}}(n)=\sum_{k=1}^{n}O(k^3)=O(n^4).
\label{oldIG}
\end{equation}

By contrast, each recurrence in \eqref{normalRec}, \eqref{gammaRec}, and \eqref{igRec} uses a fixed number of previous values. Consequently:
\begin{equation}
T_{\mathrm{rec}}(n)=\sum_{k=1}^{n}O(1)=O(n)
\qquad \text{and}\qquad
M_{\mathrm{rec}}(n)=O(1).
\label{newcomplexity}
\end{equation}

Thus, for the same sequence of $G_n^{(1)}$ approximants, the recurrence reduces the arithmetic-operation count to $O(n)$, compared with $O(n^3)$ for the Gaussian and gamma direct expressions and $O(n^4)$ for the inverse Gaussian direct expression. The bound concerns arithmetic operations; it does not include a possible increase in precision with $n$.

\section{Numerical results and discussion}
\label{sec_numerical}
The numerical experiments examine the stopping orders of the finite recurrence algorithms, their dependence on the distribution parameters, and their computational cost relative to the earlier nested-sum implementations. We use the estimated relative-error criterion used in \cite{Safouhi43}. For $n\geq2$, the estimate is:
\begin{equation}
\varepsilon_n=\left|\frac{\Gtilde_n^{(1)}(x)-\Gtilde_{n-1}^{(1)}(x)}{\Gtilde_n^{(1)}(x)}\right| ,
\label{estimatederror}
\end{equation}
with the prescribed stopping tolerance $\varepsilon=10^{-15}$. 

Tables~\ref{tab_2012gaussian}--\ref{tab_2012ig} present the recurrence-based approximations for all 35 parameter configurations considered in \cite{Safouhi43}, together with the first order $n_*$ meeting the stopping criterion \eqref{estimatederror} and the corresponding estimated relative error. The numerical calculations were performed using 50-decimal-digit arithmetic. The computed values are in excellent numerical agreement with the values reported in \cite{Safouhi43}. The algorithms evaluate the same $G_n^{(1)}$ transformation as the earlier explicit expressions, with substantially fewer arithmetic operations.

\begin{table}[!htbp]
\centering
\caption{Gaussian upper-tail approximations obtained using the first order $n_*$ satisfying the successive-approximation stopping criterion $\varepsilon_{n_*}\leq10^{-15}$.}
\label{tab_2012gaussian}
\small
\setlength{\tabcolsep}{6pt}
\begin{tabular}{@{}lrll@{}}
\toprule
$(x,\mu,\sigma)$ & $n_*$ & $\Gtilde_{n_*}^{(1)}$ & $\varepsilon_{n_*}$\\
\midrule
$(1.2,\,0,\,1)$ & 48 & $0.1150696702217084$ & $7.66\times10^{-16}$\\
$(1.6,\,0,\,1)$ & 33 & $0.05479929169955816$ & $3.16\times10^{-16}$\\
$(2,\,0,\,1)$ & 27 & $0.02275013194817922$ & $8.15\times10^{-16}$\\
$(3,\,0,\,1)$ & 19 & $0.001349898031630094$ & $7.27\times10^{-16}$\\
$(6,\,0,\,1)$ & 11 & $9.865876450376981\times10^{-10}$ & $3.68\times10^{-17}$\\
$(10,\,0,\,1)$ & 7 & $7.619853024160527\times10^{-24}$ & $6.35\times10^{-16}$\\
$(12,\,0,\,1)$ & 7 & $1.776482112077679\times10^{-33}$ & $3.07\times10^{-17}$\\
$(45,\,18,\,6)$ & 13 & $3.39767312473006\times10^{-6}$ & $2.22\times10^{-16}$\\
$(54.2,\,2,\,25)$ & 26 & $0.01839891734185767$ & $8.8\times10^{-16}$\\
$(0.3,\,0,\,1)$ & 184 & $0.3820885778110471$ & $9.99\times10^{-16}$\\
\bottomrule
\end{tabular}
\end{table}

\begin{table}[!htbp]
\centering
\caption{Gamma upper-tail approximations obtained using the first order $n_*$ satisfying the successive-approximation stopping criterion $\varepsilon_{n_*}\leq10^{-15}$.}
\label{tab_2012gamma}
\small
\setlength{\tabcolsep}{6pt}
\begin{tabular}{@{}lrll@{}}
\toprule
$(x,a,b)$ & $n_*$ & $\Gtilde_{n_*}^{(1)}$ & $\varepsilon_{n_*}$\\
\midrule
$(13,\,7,\,2)$ & 8 & $0.5265236225179999$ & $0$\\
$(15,\,7,\,2)$ & 8 & $0.3781546943234693$ & $0$\\
$(20,\,7,\,2)$ & 8 & $0.130141420882483$ & $0$\\
$(35,\,7,\,2)$ & 8 & $0.001470019774876196$ & $0$\\
$(40,\,7,\,2)$ & 8 & $0.0002551224958563007$ & $0$\\
$(45,\,7,\,2)$ & 8 & $4.079355717745713\times10^{-5}$ & $0$\\
$(50,\,7,\,2)$ & 8 & $6.106294461927904\times10^{-6}$ & $0$\\
$(60,\,7,\,2)$ & 8 & $1.173194200234696\times10^{-7}$ & $0$\\
$(120,\,7,\,2)$ & 7 & $6.292241332308505\times10^{-19}$ & $2.37\times10^{-16}$\\
$(12,\,2,\,3)$ & 3 & $0.09157819444367091$ & $0$\\
$(25.5,\,4.43,\,2.023)$ & 11 & $0.002517471973717717$ & $3.22\times10^{-16}$\\
$(45,\,5.432,\,4.5432)$ & 13 & $0.04539309469207845$ & $3.13\times10^{-16}$\\
$(14,\,1.111,\,9)$ & 49 & $0.2458730883485207$ & $8.42\times10^{-16}$\\
\bottomrule
\end{tabular}
\end{table}

\begin{table}[!htbp]
\centering
\caption{Inverse Gaussian upper-tail approximations obtained using the first order $n_*$ satisfying the successive-approximation stopping criterion $\varepsilon_{n_*}\leq10^{-15}$.}
\label{tab_2012ig}
\small
\setlength{\tabcolsep}{6pt}
\begin{tabular}{@{}lrll@{}}
\toprule
$(x,\mu,\lambda)$ & $n_*$ & $\Gtilde_{n_*}^{(1)}$ & $\varepsilon_{n_*}$\\
\midrule
$(1.5,\,1,\,1)$ & 109 & $0.1892320070000199$ & $9.87\times10^{-16}$\\
$(2,\,1,\,1)$ & 86 & $0.1145245740139932$ & $8.22\times10^{-16}$\\
$(3,\,1,\,1)$ & 61 & $0.04681207925721153$ & $8.72\times10^{-16}$\\
$(4.5,\,1,\,1)$ & 44 & $0.01430118294609312$ & $8.21\times10^{-16}$\\
$(6,\,1,\,1)$ & 35 & $0.004849882133702174$ & $8.73\times10^{-16}$\\
$(10,\,1,\,1)$ & 24 & $0.0003504145372088189$ & $8.27\times10^{-16}$\\
$(16,\,1,\,1)$ & 18 & $9.439168634947232\times10^{-6}$ & $3.88\times10^{-16}$\\
$(32,\,1,\,1)$ & 12 & $1.220065663759753\times10^{-9}$ & $3.45\times10^{-16}$\\
$(24,\,2,\,4)$ & 14 & $5.104291004380161\times10^{-7}$ & $5.78\times10^{-16}$\\
$(33.46,\,4.54,\,2.78)$ & 44 & $0.006219750083881435$ & $6.93\times10^{-16}$\\
$(23,\,6.54,\,6)$ & 57 & $0.03336361646073701$ & $9.81\times10^{-16}$\\
$(0.5,\,1,\,1)$ & 288 & $0.6350244518270296$ & $9.98\times10^{-16}$\\
\bottomrule
\end{tabular}
\end{table}

For integer shape $a$, the gamma approximants become constant after a finite transformation order, and the successive differences consequently vanish. For the cases with $a=7$ in Table~\ref{tab_2012gamma}, the first zero successive difference is at $n_*=8$; for $a=2$, it is at $n_*=3$. These zero entries are not the result of numerical underflow. The parameter configurations with noninteger shape, and the Gaussian and inverse Gaussian examples, exhibit nonterminating sequences whose estimated errors vary with the transformation order.

\begin{figure}[H]
\centering
\includegraphics[width=\linewidth,height=0.3\textheight,keepaspectratio]{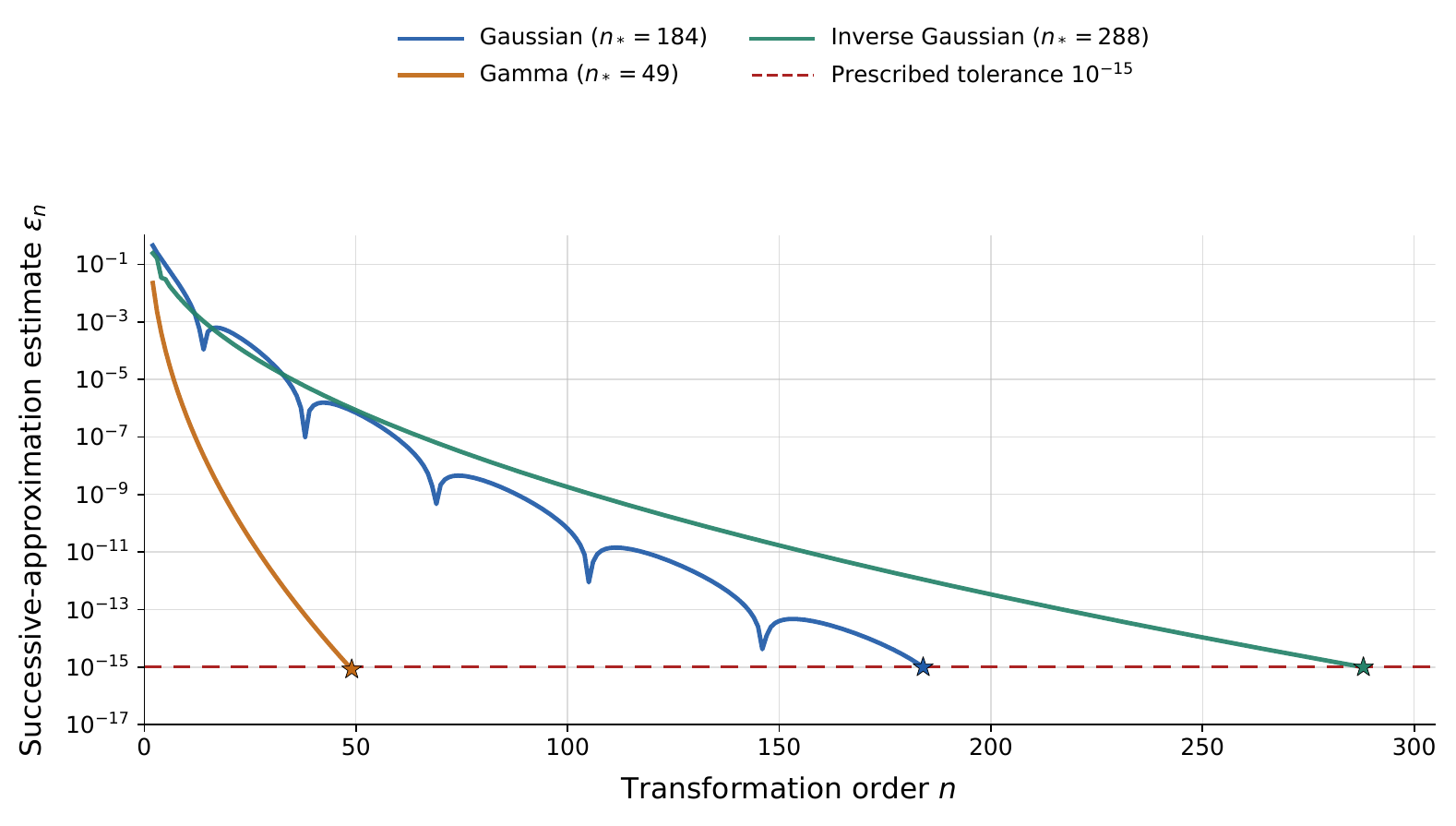}
\caption{Successive-approximation estimate $\varepsilon_n$ versus transformation order $n$ for a high-order case of each distribution: Gaussian ($x=0.3$, $\mu=0$, $\sigma=1$), gamma ($x=14$, $a=1.111$, $b=9$), and inverse Gaussian ($x=0.5$, $\mu=\lambda=1$). The examples illustrate demanding cases rather than distribution-wide behavior.}
\label{fig:conv}
\end{figure}

Figure~\ref{fig:conv} illustrates the decrease of the successive-approximation estimate $\varepsilon_n$ for one computationally demanding case of each distribution, selected from Tables~\ref{tab_2012gaussian}--\ref{tab_2012ig}. The three cases first satisfy the prescribed criterion at orders $184$, $49$, and $288$, respectively. This plot highlights the different numbers of recurrence steps required in these selected regimes, not a general ordering of the convergence rates across distributions. 

\begin{figure}[H]
\centering
\includegraphics[width=\linewidth,height=0.44\textheight,keepaspectratio]{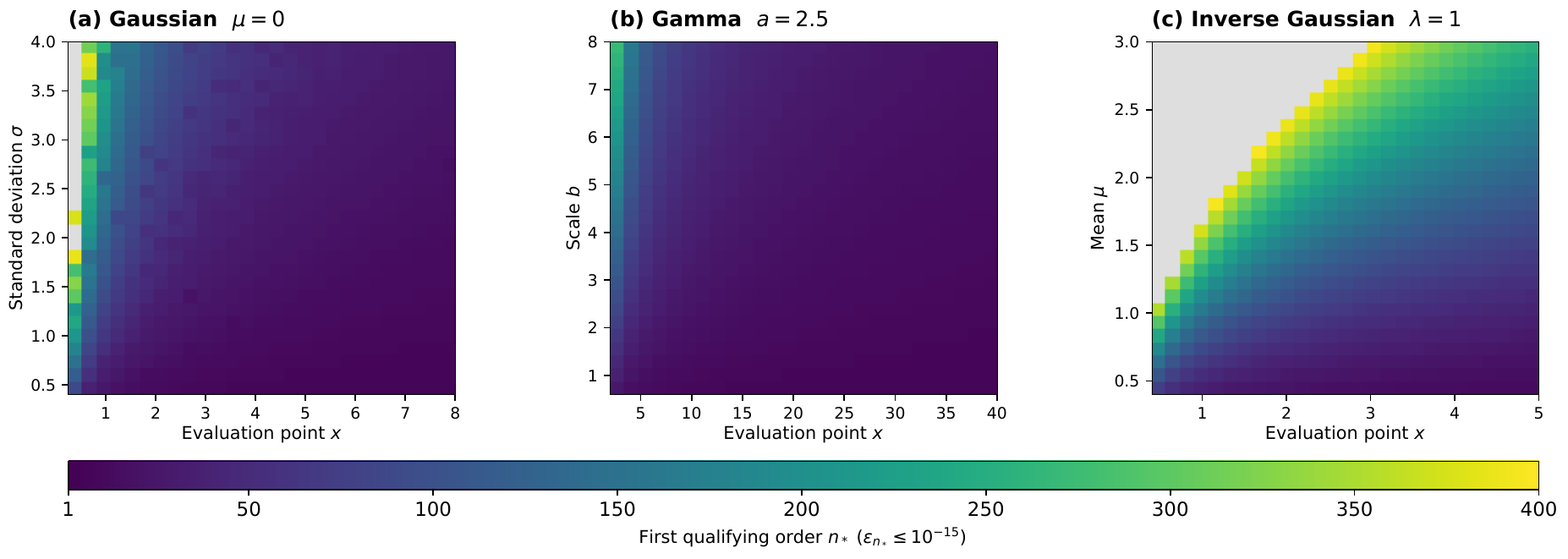}
\caption{First transformation order $n_*$ satisfying the successive-approximation tolerance $\varepsilon_{n_*}\leq10^{-15}$ on three $27\times27$ parameter grids: (a) Gaussian with $\mu=0$, (b) gamma with $a=2.5$, and (c) inverse Gaussian with $\lambda=1$. The same colour scale is used in every panel. Grey means that the criterion was not satisfied for any $n\leq400$.}
\label{fig:requiredorders}
\end{figure}

Figure~\ref{fig:requiredorders} examines the dependence of $n_*$ on two-dimensional parameter grids. The Gaussian map fixes $\mu=0$ and covers $0.25\leq x\leq8$, $0.4\leq\sigma\leq4$; the gamma map fixes $a=2.5$ and covers $2\leq x\leq40$, $0.6\leq b\leq8$; and the inverse Gaussian map fixes $\lambda=1$ and covers $0.4\leq x\leq5$, $0.4\leq\mu\leq3$. Each map contains a $27\times27$ parameter grid. The colour scale is common to all panels and records the first qualifying order when it is at most $400$. Grey cells indicate parameter combinations for which the stopping condition was not met through order $400$. The gamma panel typically requires relatively few orders over the selected parameter domain, whereas parts of the inverse Gaussian domain demand many iterations. 

\begin{figure}[H]
\centering
\includegraphics[width=\linewidth,height=0.3\textheight,keepaspectratio]{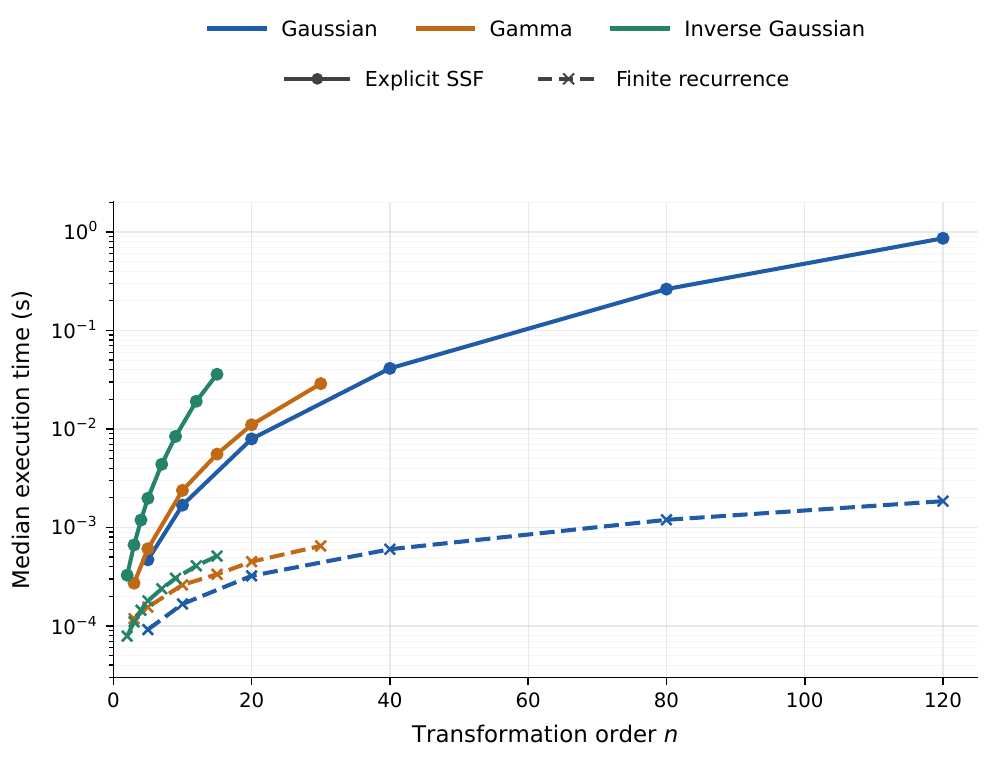}
\caption{Median execution time for generating the complete approximation sequence through order $n$ using the explicit Slevinsky--Safouhi nested-sum formulas (solid lines and circles) and the finite recurrences (dashed lines and crosses) for Gaussian, gamma, and inverse Gaussian tail probabilities. The time axis is logarithmic.}
\label{fig:timing}
\end{figure}

Figure~\ref{fig:timing} compares the execution times of the explicit formulations and the proposed recurrence algorithms. The increasing separation between the corresponding curves is consistent with the reduction in arithmetic complexity established in Section~\ref{sec_complexity}.

\FloatBarrier

\section{Conclusion}

We have derived finite recurrence relations for the $G_n^{(1)}$ transformation applied to Gaussian, gamma, and inverse Gaussian tail probabilities. The recurrences follow from a common differential structure, with an exact cancellation of the additional term arising in the reduced numerator. Consequently, the numerator and denominator satisfy recurrence relations with identical coefficients.

The resulting algorithms preserve the approximants of the original transformation while reducing the computational complexity from $O(n^3)$ and $O(n^4)$ to $O(n)$. Furthermore, only a fixed number of preceding quantities is required at each recurrence step, resulting in $O(1)$ working storage.

These improvements are particularly useful when tail probabilities must be evaluated repeatedly or at high transformation orders. Potential applications include Gaussian tail probabilities in statistical hypothesis testing, gamma tail probabilities in reliability and lifetime analysis, and inverse Gaussian tail probabilities in first-passage-time problems. In such applications, the proposed recurrences provide an efficient alternative to the explicit $G_n^{(1)}$ formulations, avoiding repeated higher-order differentiation and nested summations.

\section*{Statements and Declarations}
\noindent\textbf{Funding.} This research was supported by the Natural Sciences and Engineering Research Council of Canada (NSERC), Grant RGPIN-2016-04317 (extended to 2026).

\noindent\textbf{Competing Interests.} The authors declare that they have no competing interests.

\noindent\textbf{Data Availability.} Data sharing not applicable to this article as no datasets were generated or analyzed during the current study.


\end{document}